\documentclass[9pt]{amsart}
\usepackage{amsfonts}
\usepackage{graphicx}
\usepackage{amscd}
\usepackage{amsmath,amsfonts,amssymb,amsthm,mathrsfs,xypic}
\newtheorem{thm}{Theorem}[section]

\newtheorem{cor}[thm]{Corollary}
\newtheorem{lem}[thm]{Lemma}
\newtheorem{exm}[thm]{Example}
\newtheorem{rek}[thm]{Remark}

\newtheorem{prop}[thm]{Proposition}

\newtheorem{Def}[thm]{Definition}

\numberwithin{equation}{section}
\theoremstyle{definition}

\newcommand{\X}{\mathfrak{X}}

\begin{document}
\title [$\mathfrak{X}$-Gorenstein projective dimensions]
{$\mathfrak{X}$-Gorenstein projective dimensions}
\author [Zhibing  Zhao$^1$, Xiaowei Xu$^2$]
{Zhibing Zhao$^1$, Xiaowei Xu$^{2}$}
\thanks{The first author was partially supported by Natural Science Foundation of China (No. 11571329), the Natural Science Foundation of Anhui Province (No. 1708085MA01) and Project of University Natural Science Research of Anhui Province (No. KJ2015A101). The second author was partially supported by  Natural Science Foundation of China
(No. 11771176)}
\thanks {Email: $^1$zbzhao$\symbol{64}$ahu.edu.cn, $^2$xuxw$@$jlu.edu.cn}\maketitle
\begin{center}
$^1$School of Mathematical Sciences, Anhui University, Hefei, Anhui 230039\\
$^2$School of Mathematical Sciences, Jilin University, Changchun, Jilin 130012

\end{center}
\begin{abstract} \ In this paper, we mainly investigate the $\X$-Gorenstein projective dimension of modules and the (left) $\X$-Gorenstein global dimension of rings. Some properties of
$\X$-Gorenstein projective dimensions are obtained. Furthermore, we prove that the (left) $\X$-Gorenstein global dimension of a ring $R$ is equal to the supremum of the set of $\X$-Gorenstein
projective dimensions of all cyclic (left) $R$-modules. This result extends the well-known Auslander's theorem on the global dimension and its Gorenstein homological version.
 \vskip10pt

 \noindent 2010 Mathematics Subject Classfication: 16D80, 16E10, 16G50

\noindent Key words: $\X$-Gorenstein projective modules, $\X$-Gorenstein projective dimensions of modules, $\X$-Gorenstein global projective dimensions of rings
\end{abstract}
\section{Introduction}
Throughout this paper, $R$ always denotes a non-trivial associative ring with identity. All modules are unital left $R$-modules.
We denote the category of all left $R$-modules  by $R$-Mod. The full subcategory of $R$-Mod consisting of all Gorenstein projective $R$-modules (resp. all projective $R$-modules)
are denoted by $\mathcal{GP}(R)$ (resp. $\mathcal{P}(R)$). We use ${\rm pd}(M)$ and Gpd$(M)$ to denote, respectively, the projective dimension and Gorenstein projective dimension
of $M$.
We always assume that $\X$ is a class of $R$-modules which contains all projective $R$-modules.

As a generalization of finitely generated projective modules over a noetherian ring $R$, Auslander and Bridger \cite{AB} introduced in 1969 the notion of the module with G-dimension zero.
The G-dimension for finitely generated $R$-module $M$, denoted by G-dim$(M)$,  is defined via resolution with this kind of module. They proved that G-dim$(M)\leq {\rm pd}(M)$, with
equality G-dim$(M)={\rm pd}(M)$ when ${\rm pd}(M)$ is finite. Furthermore, they give a characterization of the Gorenstein local rings in terms of G-dimension in \cite{A67}. Namely, $R$ is Gorenstein if its
residue field $k$ has finite G-dimension, and only if every finitely generated $R$-module has finite G-dimension.

Over an arbitrary ring $R$, Enochs and Jenda defined in \cite{EO1} the Gorenstein projective module as follows\\
\noindent{\bf Definition.} An $R$-module $M$ is said to be Gorenstein projective if there exists a {\it completed projective resolution}
$$\mathbf{P}:=\cdots\rightarrow P_1\rightarrow P_0\rightarrow P^0\rightarrow P^1\rightarrow\cdots$$
being exact in $R$-$\rm Mod$ which is still exact after applying ${\rm Hom}_R(-,P)$ for any projective module $P$ and $M={\rm Ker}(P_0\rightarrow P^0)$.

The exact sequence $\mathbf{P}$ is called a complete projective resolution of $M$.\\

Avramov, Buchweitz, Martsinkovsky and Reiten proved that a finitely generated module over a noetherian ring is Gorenstein projective if and only if it is a module with G-dimension zero
(see \cite [Theorem 4.2.6]{CHris} and \cite [P.99]{CHris}). The Gorenstein projective modules share many nice properties of the classical projective module (see, for instance, \cite{CH1,CHris,EO2,H,LZ}). Gorenstein projective dimension of an arbitrary $R$-module, ${\rm Gpd}(M)$, is defined via resolution with Gorenstein projective modules.
This kind of homological dimension is analogous to the classical projective dimension and shares some of its principal properties (see \cite{BHW,BM,CHris,H} for more details).

In \cite{BO}, Bennis and Ouarghi introduced the notion of  $\X$-Gorenstein projective module where $\X$ is a class of $R$-modules that contains all projective $R$-modules (see Definition \ref{definition of X-GP module}). This notion uniforms some homological modules, including projective modules, Gorenstein projective modules and strongly Gorenstein flat modules
and so on (see Remark \ref{remark of some kinds of X-Gp modules}). Similar to the Gorenstein projective dimension, the $\X$-Gorenstein projective dimension is defined
via resolution with $\X$-Gorenstein projective modules.

In this paper, we will investigate some properties of $\X$-Gorenstein projective dimension of modules and $\X$-Gorenstein global projective dimensions of rings.

 In section 2, we will give some basic properties of $\X$-Gorenstein projective modules.

 Section 3 deal with the $\X$-Gorenstein projective dimension, $\X$-Gpd$(-)$. Firstly, we establish the following equality.\\

 \noindent{\bf Proposition.} For any family of $R$-modules $(M_i)_{i\in I}$,
$$\X{\mbox-{\rm Gpd}}(\bigoplus_{i\in I}M_i)={\rm sup}\{\X{\mbox -{\rm Gpd}}(M_i)\mid i\in I\}.$$

Note that this result is a generalization of \cite [Proposition 2.19] {H}. The following properties for $\X$-Gorenstein projective dimension are similar to those of classical projective dimension.\\

\noindent{\bf Theorem.} Let $0\rightarrow A \rightarrow B\rightarrow C\rightarrow 0$ be a short exact sequence in $R$-Mod. If two of the modules $A,B$ and $C$ have finite
$\X$-Gorenstein projective dimension, then so has the third.\\

Some relations among some homological global dimensions are also studied in this section. We show that  the (left) classical finitistic projective
dimension $l.{\rm FPD}(R)$, (left) finitistic Gorenstein projective dimension $l.{\rm FGPD}(R)$ and (left) finitistic $\X$-Gorenstein projective dimension, $l.\X{\mbox {\rm -FGPD}}(R)$, of an arbitrary ring $R$ are
equal. The last is defined as
$$l.\X\mbox{ -FGPD}(R)={\rm sup}\{\X\mbox {-Gpd}(M)\mid M \mbox{ is an $R$-module with }\X\mbox{ -Gpd}(M)<\infty\}.$$

\noindent{\bf Theorem.} For a ring $R$,  $l.{\rm FGPD}(R)=l.\X$-FGPD$(R)=$$l.{\rm FPD}(R)$.\\

In section 4, we study the left $\X$-Gorenstein global dimension of a ring $R$, $l.\X$-${\rm Ggldim}(R)$. We prove that the ideal-theoretic characterization of left
$\X$-Gorenstein global projective dimension of rings holds true.\\

\noindent{\bf Theorem.} For a ring $R$,
 \begin{align*}l.\X{\mbox -{\rm Ggldim}} (R)&={\rm sup}\{\mbox{$\X$-}{\rm Gpd}(M)\mid \mbox {$M$ is a finitely generated left $R$}\mbox {-module}\}\\
&={\rm sup}\{\mbox{$\X$-}{\rm Gpd}(R/I)\mid \mbox {$I$ is a left ideal of $R$}\}.\\
\end{align*}

As a consequence, we can obtain the well-known Auslander's theorem on the global dimension of rings and the Gorenstein counterpart of the Auslander's theorem.

\section{$\X$-Gorenstein projective modules}

In this section, we will review some properties of $\X$-Gorenstein projective modules. The following definition was introduced firstly by Bennis and Ouarghi in \cite{BO}.

\begin{Def}\label{definition of X-GP module}$\rm ($\cite{BO}$\rm )$ Let $\mathfrak{X}$ be a class of left $R$-modules that contains all projective left $R$-modules. A left $R$-module $G$ is called $\mathfrak{X}$-Gorenstein projective
if there exists an exact sequence
$$\mathbf{P}:=\cdots\rightarrow P_1\rightarrow P_0\rightarrow P^0\rightarrow P^1\rightarrow\cdots$$
of projective left $R$-modules which is still exact after applying ${\rm Hom}_A(-,X)$ for any module $X\in\X $ and $G={\rm Ker}(P_0\rightarrow P^0)$.
Furthermore, the exact sequence $\mathbf{P}$ is called an $\X$-complete projective resolution of $G$. We denoted by $\X$-$\mathcal{GP}(R)$ the full subcategory of $R$-$\rm Mod$ consisting
of all $\X$-Gorenstein projective modules.

\end{Def}
\begin{rek}\label{remark of some kinds of X-Gp modules}$\rm(1)$ It is easy to see that
$$\mathcal{P}(R)\subseteq \X{\rm-}\mathcal{GP}(R)\subseteq\mathcal{GP}(R).$$

$\rm(2)$ Let $\X$ be the class of projective $R$-modules. Then $\X$-$\mathcal{GP}(R)=\mathcal{GP}(R)$.

$\rm (3)$ Let $\X$ be the class of Gorenstein projective $R$-modules. Then $\X$-$\mathcal{GP}(R)=\mathcal{P}(R)$$\rm($see \cite{T}$\rm)$.

$\rm(4)$ Let $\X$ be the class of flat $R$-mdoules. Then the class of $\X$-Gorenstein projective $R$-modules coincide with the class of strongly
Gorenstein flat $R$-modules $\rm($see \cite{DLM}$\rm)$.

\end{rek}
\begin{exm}$(1)$ Put $R=k[x]/(x^2)$, where $k$ is a field. Then $R$ is a local quasi-Frobenius ring. Let $(\bar{x})$ be the residue class of $x$ in $R$ and $\X$ a
class of modules which contains all left projective $R$-modules and $(\bar{x})$. Then $(\bar{x})$ is a Gorenstein projective module but it is not $\X$-Gorenstein projective module.
$\rm($for details see \cite [Example 3.7] {T}$\rm)$

$(2)$ Put $R=\mathbb{Z}/4\mathbb{Z}$, where $\mathbb{Z}$ is the ring of integers. Then $R$ is a quasi-Frobenius ring. Let $\X$ be the class of flat $R$-modules.
Then $2R$ is an $\X$-Gorenstein projective module, but it is not projective.
\end{exm}
Using the definition of $\X$-Gorenstein projective module, we immediately get the following characterization.
\begin{prop}\label{Characterization of X-GP module}Let $G$ be a left $R$-module. Then the followings are equivalent.

$(1)$ G is $\X$-Gorenstein projective.

$(2)$ ${\rm i)}$ ${\rm Ext}^i_R(G,X)=0$ for any $X\in\X$ and $i>0$;

\quad ${\rm ii\rm)}$ There exists an exact sequence of $R$-modules $\mathbf{Q}:= 0\rightarrow G\rightarrow P^0\rightarrow P^1\rightarrow \cdots$ with $P^i$ projective for every $i\geq 0$
such that ${\rm Hom}(\mathbf{Q},X)$ is still exact for every $X\in\X$.

$(3)$ There exists a short exact sequence of $R$-modules $0\rightarrow G\rightarrow P\rightarrow G'\rightarrow 0$, where $P$ is projective and $G'$ is $\X$-Gorenstein projective.
\end{prop}

In fact, we can always assume that the modules in an $\X$-complete projective resolution are free.
\begin{prop} Let $G$ be an $\X$-Gorenstein projective module. Then there exists an $\X$-complete projective resolution $\mathbf{F}:=\cdots\rightarrow F_1\rightarrow F_0
\rightarrow F^0\rightarrow F^1\rightarrow \cdots$ with $F_i$ and $F^i$ being free such that $G\cong{\rm Im}(F_0\rightarrow F^0)$.
\end{prop}

\noindent{\bf Proof.} The proof is similar to the one of \cite [Proposition 2.4]{H}. For the sake of completeness, we give the proof as follows.

 It suffices to  construct the ``right half" of $\mathbf{F}$. By Proposition \ref{Characterization of X-GP module}, There exists an exact sequence of $R$-modules
 $\mathbf{Q}:= 0\rightarrow G\rightarrow P^0\rightarrow P^1\rightarrow \cdots$ with $P^i$ projective for every $i\geq 0$
such that ${\rm Hom}(\mathbf{Q},X)$ is still exact for every $X\in\X$. We can pick projective modules $Q^0,Q^1,Q^2,\cdots$ such that all of the modules
$F^0=P^0\oplus Q^0 $ and $F^i=P^i\oplus Q^{i-1}\oplus Q^i$ are free for $i>0$.

By adding $0\longrightarrow
Q^i\stackrel{=}\longrightarrow Q^i\longrightarrow 0$ to $\mathbf{Q}$ in degree $i$ and $i+1$, we obtain the desired sequence. \hfill$\square$\\

Recall that a class $\X$ of $R$-modules is called {\it projectively resolving} if $\mathcal{P}(R)\subseteq\X$, and for every short exact sequence $0\rightarrow A\rightarrow B
\rightarrow C\rightarrow 0$ in $R$-$\rm Mod$ with $C\in\X$ the conditions $A\in\X$ and $B\in\X$ are equivalent. The following proposition comes from \cite{BO}.
\begin{prop}\label{X-GP subcat is projective resolving}$\rm($\cite[Theorem 2.3]{BO}$\rm)$ Let $0\rightarrow A \rightarrow B\rightarrow C\rightarrow 0$ be a short exact sequence in $R$-$\rm Mod$. Then we have

$(1)$ If $A$ and $C$ are $\X$-Gorenstein projective, then so is $B$.

$(2)$ If the exact sequence is split and $B$ is $\X$-Gorenstein projective, then $A$ and $C$  are also $\X$-Gorenstein projective.

$(3)$ If $B$ and $C$ are $\X$-Gorenstein projective, then so is $A$.

Hence the class of $\X$-$\mathcal{GP}(R)$ of all $\X$-Gorenstein projective modules is projectively resolving. Furthermore, $\X$-$\mathcal{GP}(R)$ is closed under arbitrary direct sums and direct summands.
\end{prop}
\noindent{\bf Proof.} The last two assertions are easy from (1), (2) and (3). The proof is similar to the one of \cite [Theorem 2.5] {H}. Here, we give a new and simple proof as follows.

(1) Since $A$,$C$ are $\X$-Gorenstein projective modules, there exist the following exact sequences
$0\rightarrow A\rightarrow P^0\rightarrow A^1\rightarrow 0$ and
$0\rightarrow C\rightarrow Q^0\rightarrow C^1\rightarrow 0$
where $P^0, Q^0$ are projective and $A^1, C^1$ are $\X$-Gorenstein projective.

From the Long Exact Sequence Theorem and  Propostion \ref{Characterization of X-GP module}, we get ${\rm Ext}^i_R(B,X)=0$ for any $X\in\X$ and any $i>0$. Then we have the following
commutative diagram with exact rows and columns
$$\xymatrix{
         & 0 \ar[d]  & 0 \ar[d] & 0 \ar[d]  \\
      0  \ar[r]  & A \ar[d] \ar[r] & B \ar[d]^{f} \ar[r] & C \ar[d] \ar[r] & 0 \\
      0  \ar[r] &  P^0 \ar[d] \ar[r] & P^0\oplus Q^0 \ar[d] \ar[r] & Q^0 \ar[d]\ar[r] & 0  \\
       & A^1  \ar[d]  & B^1  \ar[d]  & C^1  \ar[d] \\
 & 0  & 0  & 0  }$$
 where $B^1={\rm Coker}f$. By the Snake Lemma, there is an exact sequence $0\rightarrow A^1\rightarrow B^1\rightarrow C^1\rightarrow
0$  with $A^1$ and $C^1$ being $\X$-Gorenstein projective. Then we have ${\rm Ext}^i_R(B^1,X)=0$ for any $X\in\X$ and any $i>0$ by Propostion \ref{Characterization of X-GP module}.
Repeating this process, we get an exact sequence $0\rightarrow B\rightarrow P^0\oplus Q^0\rightarrow P^1\oplus Q^1\rightarrow \cdots$ with $P^i,Q^i$ projective for every $i\geq 0$
such that which  is still exact after applying ${\rm Hom}(-,X)$ for every $X\in\X$. It follows from Proposition \ref{Characterization of X-GP module} that $B$ is $\X$-Gorenstein projective.

(2) Assume that $B\cong A\oplus C$ is $\X$-Gorenstein projective. Then ${\rm Ext}^i_R(B,X)=0$ for any $X\in\X$ and any $i>0$, and there is an exact sequence
 $\xymatrix@C=0.5cm{
  0 \ar[r] & A\oplus C \ar[r]^{f} & P^0 \ar[r] & L^1 \ar[r] & 0 }$, where $P^0$ is projective and $L^1$ is $\X$ Gorenstein projective.
  Put $A^1=\hbox{\rm Coker}(f\circ i_A)$, that is, $0\rightarrow A\rightarrow P^0\rightarrow A^1\rightarrow
  0$ is exact. Hence there is an exact sequence $0\rightarrow C\rightarrow A^1\rightarrow L^1\rightarrow
  0$. From the Long Exact Sequence Theorem, ${\rm Ext}^i_R(A^1,X)=0$ for any $X\in\X$ and any $i>0$ and $\xymatrix@C=0.5cm{
  0 \ar[r] & A\oplus C \ar[r] & A\oplus A^1 \ar[r] & L^1 \ar[r] & 0}$ is exact.

  By (1), $A\oplus A^1$ is $\X$-Gorenstein projective. Similar to the process above, we get an exact sequence $0\rightarrow A^1\rightarrow
P^1\rightarrow A^2\rightarrow 0$, where $P^1$ is projective and $A^2$ is $\X$-Gorenstein projective. Repeating this process, we will obtain a long exact sequence
$0\rightarrow A\rightarrow P^0\rightarrow P^1\rightarrow\cdots$, which is the desired long exact sequence in
Proposition \ref{Characterization of X-GP module} (2) ii). Therefore, $A$ is $\X$-Gorenstein projective. Similarly, we can prove that $C$ is also $\X$-Gorenstein projective.

(3) Since $C$ is $\X$-Gorenstein projective, there is an exact sequence $0\rightarrow
C_1\rightarrow P_0\rightarrow C\rightarrow 0$ where $P_0$ is projective and $C_1$ is $\X$-Gorentsein projective.

Consider the following pullback diagram
$$\xymatrix{
       &   & 0 \ar[d]  & 0 \ar[d]   \\
       &  & C_1 \ar[d] \ar@{=}[r] & C_1 \ar[d]  \\
      0  \ar[r] & A \ar@{=}[d] \ar[r] & E \ar[d] \ar[r] & P_0 \ar[d] \ar[r] & 0  \\
      0  \ar[r] & A  \ar[r] & B  \ar[d] \ar[r] & C \ar[d] \ar[r] & 0  \\
      & & 0  & 0   }$$
where the second row is split, that is, $E\cong A\oplus
P_0$. By (1) and the second column, $E$ is $\X$-Gorenstein projective. It follows from $A$ is a direct summand of $E$ that $A$ is $\X$-Gorenstien projective.
\hfill${\square}$\\

By the Proposition above, the following result is a direct consequence of \cite [Lemma 3.12] {AB}.

\begin{prop}\label{X-Gpd similar to classical H-dimension} Let $M$ be an $R$-module and consider two exact sequences,
$$0\rightarrow K_n\rightarrow G_{n-1}\rightarrow\cdots\rightarrow G_0\rightarrow M\rightarrow 0,$$
$$0\rightarrow K'_n\rightarrow G'_{n-1}\rightarrow\cdots\rightarrow G'_0\rightarrow M\rightarrow 0,$$
where $G_0,\cdots,G_{n-1}$ and $G'_0,\cdots,G'_{n-1}$ are $\X$-Gorenstein projective modules. Then $K_n$ is $\X$-Gorenstein projective if and only if $K'_n$ is $\X$-Gorenstein projective.

\end{prop}

\section{$\X $-Gorenstein projective dimensions}
In this section, we will investigate some properties of the $\X$-Gorenstein projective dimension, $\X$-Gpd$(-)$. It is defined as follows
\begin{Def}\label{definition of X-Gorenstein dimension} Let $M$ be an $R$-module. The $\X$-Gorenstein projective dimension of $M$, denoted by $\X$-${\rm Gpd}(M)$, is defined as
 $\X$-${\rm Gpd}(M)={\rm Inf}\{n\mid \exists \X$-Gorenstein projective resolution $0\rightarrow G_n\rightarrow \cdots G_1\rightarrow G_0\rightarrow M\rightarrow 0  \}$. We set
 $\X$-${\rm Gpd}(M)=\infty$ if there is no such an $\X$-Gorenstein projective resolution.

 We define $l.\X$-${\rm Ggldim}(R)={\rm sup}\{\X$-${\rm Gpd}(M)\mid M$ is any left $R$-module $\}$, and call it left  $\X$-Gorenstein global projective dimension of $R$.
\end{Def}
\begin{rek} $(1)$ It is obvious that ${\rm Gpd}(M)\leq \X$-${\rm Gpd}(M)\leq {\rm pd}(M)$ for any left $R$-module $M$.

$(2)$ If $\X$-${\rm Gpd}(M)<\infty$, then $\X$-${\rm Gpd}(M)={\rm Gpd}(M)$ $($See \cite [Corollary 3.5] {MP}$)$.
\end{rek}
Let $\X$ be a class of $R$-modules which contains all projective $R$-modules. For any $M\in R$-Mod, there exists an exact sequence
$$\cdots \rightarrow X_n\rightarrow\cdots\rightarrow X_1\rightarrow X_0\rightarrow M\rightarrow0$$
with $X_n\in\X$ for $n=0,1,\cdots$. Similar to Definition \ref{definition of X-Gorenstein dimension}, we can define $\X$-dimension of $M$.
We denote by $\widetilde{\X}$ the class of all $R$-module which has finite $\X$-dimension. In particular, $\widetilde{\mathcal{P}}(R)$ (resp. $\widetilde{\mathcal{GP}}(R)$)
denotes the class of the $R$-module which has finite projective dimension (resp. the class of the $R$-module which has finite Gorenstein projective dimension).
The following proposition comes from \cite{MP}, which is a characterization of $\X$-${\rm Gpd}(M)$.
\begin{prop}\label{characterization of X-Gpd} Let $M$ be a left $R$-module with finite $\X$-Gorenstein projective dimension. Then the following statements are equivalent for a nonnegative integer $n$.

$(1)$ $\X$-${\rm Gpd}(M)\leq n$.

$(2)$ ${\rm Ext}^i_R(M,X)=0$ for all $i>n$, and any $R$-modules $X$ in $\X$.

$(3)$ ${\rm Ext}^i_R(M,\widetilde{X})=0$ for all $i>n$, and any $R$-modules $\widetilde{X}\in \widetilde{\X}$.

$(4)$ For every exact sequence $0\rightarrow K_n\rightarrow G_{n-1}\rightarrow\cdots\rightarrow G_0\rightarrow M\rightarrow 0$ with $G_0,G_1,\cdots G_{n-1}$ are
$\X$-Gorenstein projective, then also $K_n$ is $\X$-Gorenstein projective.
\end{prop}
\begin{prop}If $(M_i)_{i\in I}$ is a family of $R$-modules, then \\

\quad\quad\quad\quad\quad\quad\quad\quad\quad $\X$-${\rm Gpd}(\bigoplus_{i\in I}M_i)={\rm sup}\{\X$-${\rm Gpd}(M_i)\mid i\in I\}$.
\end{prop}
\noindent{\bf Proof.} The inequality $\X$-Gpd$(\bigoplus_{i\in I}M_i)\leq{\rm sup}\{\X$-${\rm Gpd}(M_i)\mid i\in I\}$ is clear since $\X$-$\mathcal{GP}(R)$ is closed under direct sums.

For the converse inequality ``$\geq$", it suffices to show that $\X$-Gpd$(M')\leq\X$-Gpd$(M)$ for any direct summand $M'$ of $M$.
Without loss of generality, we assume that $\X$-Gpd$(M)=n$ is finite. We do by induction on $n$.

If $n=0$,  $M$ is $\X$-Gorenstein projective, then $M'$ is also $\X$-Gorenstein projective by Proposition \ref{X-GP subcat is projective resolving}.

Now, we assume that $n>0$. We write that $M=M'\oplus M''$ for some $R$-module $M''$. Pick exact sequences $0\rightarrow K'\rightarrow P'\rightarrow M'\rightarrow0$
and $0\rightarrow K''\rightarrow P''\rightarrow M''\rightarrow0$, where $P'$ and $P''$ are projective. By the Horseshoe Lemma, we get the following commutative diagram
$$\xymatrix{
  & 0 \ar[d] & 0 \ar[d]  & 0 \ar[d]  \\
  0  \ar[r] & K' \ar[d] \ar[r] & K'\oplus K'' \ar[d] \ar[r] & K'' \ar[d] \ar[r] & 0 \\
  0  \ar[r] & P' \ar[d] \ar[r] & P'\oplus P'' \ar[d] \ar[r] & P'' \ar[d] \ar[r] & 0 \\
  0  \ar[r] & M' \ar[d] \ar[r] & M \ar[d] \ar[r] & M'' \ar[d] \ar[r] & 0  \\
  & 0  & 0  & 0  }$$
with exact columns and split exact rows. It follows from Proposition \ref{X-Gpd similar to classical H-dimension} that $\X$-Gpd$(K'\oplus K'')=n-1$. Hence the induction hypothesis yields that $\X$-Gpd$(K')\leq n-1$. Therefore, $\X$-Gpd$(M')\leq n$ by
Proposition \ref{X-Gpd similar to classical H-dimension} again. We complete the proof.   \hfill$\square$

Similar to the case of the projective dimension of modules, we will show that if two terms of a short exact sequence  have finite
$\X$-Gorenstein projective dimension, then so has the third.

\begin{lem}\label{syzygy is X-GP} Let $M$ be an $R$-module. If $M$ has an $\X$-Gorenstein projective resolution $\mathbf{G}\rightarrow M$ of length at most $n$, then in each projective resolution $\mathbf{P}\rightarrow M$, the $i$th syzygy module $\Omega^i(\mathbf{P})$ is $\X$-Gorenstein projective for all $i\geq n$.
\end{lem}
\noindent{\bf Proof.} It is easy to see that there is a comparison map $\gamma:\mathbf{P}\rightarrow\mathbf{G}$ being a quasi-isomorphism. Hence its mapping cone is exact.
Therefore, we have an exact sequence of $R$-modules
$$0\rightarrow \Omega^i(\mathbf{P})\rightarrow P_{i-1}\rightarrow\cdots\rightarrow P_n\rightarrow P_{n-1}\oplus G_n\rightarrow\cdots P_0\oplus G_1\rightarrow G_0\rightarrow 0.$$
By Proposition \ref{X-GP subcat is projective resolving}, $P_{j-1}\oplus G_j$ is $\X$-Gorenstein projective for each $j=1,2,\cdots,n$. It follows from Proposition \ref{X-GP subcat is projective resolving} again that $\Omega^i(\mathbf{P})$ is $\X$-Gorenstein projective for all $i\geq n$.   \hfill$\square$

\begin{thm}\label{finity of X-Gpd of short exact sequence}Let $0\rightarrow A \rightarrow B\rightarrow C\rightarrow 0$ be a short exact sequence in $R$-$\rm Mod$. We have

$(1)$ If $A$ and $C$  have finite $\X$-Gorenstein projective dimension, then so has $B$.

$(2)$ If $A$ and $B$  have finite $\X$-Gorenstein projective dimension, then so has $C$.

$(3)$ If $B$ and $C$  have finite $\X$-Gorenstein projective dimension, then so has $A$.
\end{thm}
\noindent{\bf Proof.} $(1)$ Pick projective resolutions $\mathbf{P}_A\rightarrow A$ and $\mathbf{P}_C\rightarrow C$ of $A$ and $C$, respectively. By the Horseshoe Lemma,
$\mathbf{P}_B=\mathbf{P}_A\oplus\mathbf{P}_C$ is a projective resolution of $B$. The split exact sequence of complexes
$0\rightarrow \mathbf{P}_A\rightarrow \mathbf{P}_B\rightarrow \mathbf{P}_C\rightarrow 0$ induces, for every $n\in\mathbb{N}$, an exact sequence of $R$-modules
$$0\rightarrow \Omega^n(\mathbf{P}_A)\rightarrow\Omega^n(\mathbf{P}_B)\rightarrow\Omega^n(\mathbf{P}_C)\rightarrow 0.$$
Since $\X$-$\mathcal{GP}(R)$ is closed under extensions, Applying Lemma \ref{syzygy is X-GP}, we get that $\Omega^i(\mathbf{P}_B)$ is $\X$-Gorenstein projective when
$i\geq{\rm max}\{\X{\mbox-{\rm Gpd}}(A), \X{\mbox-{\rm Gpd}}(C) \}$. Hence
$\X$-Gpd$(B)$ is finite.

 $(2)$ Assume that  $\X$-Gpd$(A)$ and $\X$-Gpd$(B)$ are finite. There exists an exact sequence
 $0\rightarrow K\rightarrow G\rightarrow B\rightarrow 0$ with $G$ $\X$-Gorenstein projective and $\X$-Gpd$(K)<\infty$.

 Consider the following pullback diagram
$$\xymatrix{
 & 0 \ar[d]  & 0 \ar[d] \\
  & K \ar[d] \ar@{=}[r] & K \ar[d]  \\
  0 \ar[r] & T \ar[d] \ar[r] & G \ar[d] \ar[r] & C \ar@{=}[d] \ar[r] & 0  \\
  0  \ar[r] & A \ar[d] \ar[r] & B \ar[d] \ar[r] & C  \ar[r] & 0  \\
  & 0  & 0    }$$
since $\X$-Gpd$(K)$ and $\X$-Gpd$(A)$ are finite, $\X$-Gpd$(T)$ is also finite. It follows from the second exact row that $\X$-Gpd$(C)$ is finite.

$(3)$ If $\X$-Gpd$(B)$ is finite, then there exists an exact sequence
$0\rightarrow K\rightarrow G\rightarrow B\rightarrow 0$ with $G$ $\X$-Gorenstein projective and $\X$-Gpd$(K)<\infty$. Then we can get the same pullback diagram
above.
since $\X$-Gpd$(C)$ is finite, $\X$-Gpd$(T)$ is also finite. Considering the first exact column in the diagram above, we get $\X$-Gpd$(A)<\infty$ by $(2)$. \hfill$\square$

\begin{cor}\label{X-Pd of short exact sequence}Let $0\rightarrow A \rightarrow B\rightarrow C\rightarrow 0$ be a short exact sequence in $R$-${\rm Mod}$. We have

$(1)$ $\X$-${\rm Gpd}(A)\leq{\rm sup}\{\X$-${\rm Gpd}(B)$, $\X$-${\rm Gpd}(C)-1\}$ with equality if $\X$-${\rm Gpd}(B)\neq\X$-${\rm Gpd}(C)$.

$(2)$ $\X$-${\rm Gpd}(B)\leq{\rm sup}\{\X$-${\rm Gpd}(A)$, $\X$-${\rm Gpd}(C)\}$ with equality if $\X$-${\rm Gpd}(C)\neq\X$-${\rm Gpd}(A)+1$.

$(3)$ $\X$-${\rm Gpd}(C)\leq{\rm sup}\{\X$-${\rm Gpd}(B)$, $\X$-${\rm Gpd}(A)+1\}$ with equality if $\X$-${\rm Gpd}(A)\neq\X$-${\rm Gpd}(B)$.
\end{cor}
\noindent{\bf Proof.} Using Proposition \ref{characterization of X-Gpd} and Theorem \ref{finity of X-Gpd of short exact sequence}, the argument is analogous to the case
of projective dimension in classical homological algebra.    \hfill$\square$

\begin{prop}\label{precover of a finite X-Gpd module}Let $M$ be an $R$-module with finite $\X$-Gorenstein projective dimension $n(\geq1)$. Then there exists an exact sequence $0\rightarrow K\rightarrow G\rightarrow M\rightarrow 0$, where $G$ is $\X$-Gorenstein projective and ${\rm pd}(K)=n-1$.
\end{prop}
\noindent{\bf Proof.} Since $\X$-Gpd$(M)=n$, there exists an exact sequence
$$0\rightarrow L\rightarrow P_{n-1}\rightarrow \cdots\rightarrow P_0\rightarrow M\rightarrow 0,$$
where $L$ is $\X$-Gorenstein projective and $P_i$ is projective for each $i=0,1,\cdots,n-1$. By the definition of $\X$-Gorenstein projective module, there exists
an exact sequence $$0\rightarrow L\rightarrow P^0\rightarrow \cdots\rightarrow P^{n-1}\rightarrow G\rightarrow 0,$$
where $P^0,\cdots,P^{n-1}$ are projective, $G$ is $\X$-Gorenstein projective, and such that the functor ${\rm Hom}(-,X)$ leaves this sequence still exact for any $X\in\X$.
Similar to the Compare Lemma, we get a chain map between complexes
$$\xymatrix{
  0  \ar[r] & P^0 \ar[d] \ar[r] & \cdots  \ar[r] & P^{n-1} \ar[d] \ar[r] & G \ar[d] \ar[r] & 0 \\
  0 \ar[r] & P_{n-1} \ar[r] & \cdots \ar[r] & P_0 \ar[r] & M \ar[r] & 0   }$$
which is a quasi-isomorphism. Its mapping cone is exact, namely,
$0\rightarrow P^0\rightarrow P_{n-1}\oplus P^1\rightarrow\cdots\rightarrow P_0\oplus G\rightarrow M\rightarrow 0$ is exact. Hence we get an short exact sequence
$0\rightarrow K\rightarrow P_0\oplus G\rightarrow M\rightarrow 0$, where $K={\rm Ker}(P_0\oplus G\rightarrow M)$ and ${\rm pd}(K)\leq n-1$, as desired. \hfill$\square$

\begin{prop}\label{finite dimension equality} If an $R$-module $M$ has finite projective dimension, then ${\rm Gpd}(M)=\X$-${\rm Gpd}(M)={\rm pd}(M)$. In particular,  $\X$-$\mathcal{GP}(R)\cap\widetilde{\mathcal{P}}(R)=\mathcal{GP}(R)\cap\widetilde{\mathcal{P}}(R)=\mathcal{P}(R)$ and
 $\mathcal{GP}(R)\cap\X$-$\widetilde{\mathcal{GP}}(R)=\X$-$\mathcal{GP}(R)$.
\end{prop}
\noindent{\bf Proof.} By  \cite [Proposition 2.27]{H},  Gpd$(M)=$pd$(M)$ for an $R$-module $M$ with finite projective dimension. There is always an inequality
$\X$-Gpd$M\leq {\rm pd}(M)$, and consequently, we also have $\X$-Gpd$(M)<\infty$. It follows from  \cite [Corollary 3.15]{MP}, that Gpd$(M)=\X$-Gpd$(M)$.
\hfill$\square$\\

We end this section with some relations among homological dimensions. The {\it left finitistic $\X$-Gorenstein projective dimension} of a ring $R$ is defined as follows

$l.\X$-${\rm FGPD}(R)={\rm sup}\{\X$-Gpd$(M)\mid M$ is an $R$-module with $\X$-Gpd$(M)<\infty\}$.

We use $l.{\rm FPD}(R)$ and $l.{\rm FGPD}(R)$ to denote the usual {\it left finitistic projective dimemsion of $R$} and {\it left finitistic Gorenstein projective dimension of $R$}, respectively.
Now, we have the following theorem.
\begin{thm} For a ring $R$, l.{\rm FGPD}$(R)=$l.$\X$-{\rm FGPD}$(R)=$l.{\rm FPD}$(R)$.
\end{thm}
\noindent{\bf Proof.} It follows from \cite [Theorem 2.28]{H} that $l.{\rm FGPD}(R)=l.{\rm FPD}(R)$. We only need to show that $l.{\rm FPD}(R)=l.\X$-${\rm FGPD}(R)$.

Clearly $l.{\rm FPD}(R)\leq l.\X$-${\rm FGPD}(R)$ by Proposition \ref{finite dimension equality}. For any $M\in R$-Mod with finite $\X$-Gorenstein projective dimension $n$, by Proposition \ref{precover of a finite X-Gpd module}, there exists a short exact sequence $0\rightarrow K\rightarrow G\rightarrow M\rightarrow 0$, where $G$ is $\X$-Gorenstein projective and ${\rm pd}(K)=n-1$. Then $l.\X$-${\rm FGPD}(R)\leq l.{\rm FPD}(R)+1$.

 Now, we can assume that $l.\X$-FGPD$(R)=n<\infty$, and then pick an $R$-module $M$
such that $\X$-Gpd$(M)=n$. By Proposition \ref{precover of a finite X-Gpd module}, there exists an exact sequence $0\rightarrow K\rightarrow G\rightarrow M\rightarrow0$, where
$G$ is $\X$-Gorenstein projective and ${\rm pd}(K)=n-1$.

If $G$ is projective, then ${\rm pd}(M)=n$, and hence $l.{\rm FPD}(R)=n=l.\X$-FGPD$(R)$.

 If $G$ is not a projective module, there exists a projective module $P$ such that $G\subsetneqq P$. Put $L=P/K$, there is an exact sequence
 $$0\rightarrow M(\cong G/K)\rightarrow L\rightarrow L/M(\cong P/G)\rightarrow 0.$$
 If $L$ is $\X$-Gorenstein projective, then $\X$-Gpd$(L/M)=n+1$, this contradict to the fact that $l.\X$-FGPD$(R)=n<\infty$. Hence $L$ is not $\X$-Gorenstein projective.
 In particular, $L$ is not projective. It follows from the short exact sequence $0\rightarrow K\rightarrow P\rightarrow L\rightarrow0$ that
 ${\rm pd}(L)={\rm pd}(K)+1=n$. Therefore we get that $l.{\rm FPD}(R)=n=l.\X$-FGPD$(R)$.\hfill$\square$

\section{$\X$-Gorenstein global dimension of rings}

In this section, the (left) $\X$-Gorenstein global projective dimension of rings is studied. We prove that the (left) $\X$-Gorenstein global dimension of ring $R$ is equal to
 the supremum of the set of $\X$-Gorenstein projective dimensions of all cyclic $R$-modules. In order to prove the main theorem, we need some lemmas

\begin{lem}\label{another characterization of X-Gpd} Let $M$ be an $R$-module and $n$ a non-negative integer. Then, $\X$-${\rm Gpd}(M)\leq n$ if and only if there exists a short exact sequence
of $R$-modules, $0\rightarrow M\rightarrow H\rightarrow G\rightarrow 0$, where ${\rm pd}(H)\leq n$ and $G$ is $\X$-Gorenstein projective.
\end{lem}
\noindent{\bf Proof.} The sufficiency follows from Corollary \ref{X-Pd of short exact sequence} (1).

Now, we suppose that $\X$-${\rm Gpd}(M)\leq n$. By Proposition \ref{precover of a finite X-Gpd module}, there exists a short exact sequence $0\rightarrow K\rightarrow G'\rightarrow M\rightarrow 0$, where $G'$ is $\X$-Gorenstein projective and ${\rm pd}(K)=n-1$. Since $G'$ is $\X$-Gorenstein projective, there exists a short exact sequence
$0\rightarrow G'\rightarrow P\rightarrow G\rightarrow 0$, where $P$ is projective and $G$ is $\X$-Gorenstein projective. Hence we get the following pushout diagram
$$\xymatrix{
       &   & 0 \ar[d]  & 0 \ar[d]   \\
        0  \ar[r] &  K \ar@{=}[d] \ar[r] & G' \ar[d] \ar[r] & M \ar[d] \ar[r] & 0 \\
      0  \ar[r] & K  \ar[r] & P \ar[d] \ar[r] & H \ar[d] \ar[r] & 0  \\
       &    & G  \ar[d] \ar@{=}[r] & G \ar[d]  \\
      & & 0  & 0   }$$
where ${\rm pd}(H)={\rm pd}(K)+1\leq n$. Therefore, the third exact column is the exact sequence we desired. \hfill$\square$

\begin{lem}\label{generalization of Baer criterion} Let $n$ be a nonnegative integer. If $\X\mbox{\rm -Gpd}(R/I)\leq n$ for any left ideal $I$ of $R$, then ${\rm id}(X)\leq n$
 for $X\in\X$, where ${\rm id}(X)$ denote the injective dimension of $X$.
\end{lem}
\noindent{\bf Proof.} It follows from Baer's Criterion of injective dimension \cite [Theorem 3.30] {R}. \hfill$\square$

The following Eklof's Lemma is well-known (see \cite{Ek} or \cite [Lemma 2.5] {FS}).
\begin{lem}\label{Eklof Lemma} Suppose that, for an ordinal $\tau$,
$$0=M_0 <M_1<\cdots<M_{\alpha}<M_{\alpha+1}<\cdots (\alpha<\tau)$$
is well-ordered continuous ascending chain of submodules of an $R$-module $M$ whose union is $M$. If, for some integer $i>1$
and some $R$-module $X$, ${\rm Ext}^i_R(M_{\alpha+1}/M_{\alpha},X)=0$ for all $\alpha+1<\tau$, then ${\rm Ext}^i_R(M,X)=0$.
\end{lem}

\begin{thm}\label{X-Ggldim of R} Let $R$ be a ring and $n$ a non-negative integer. The following statements are equivalent.

$(1)$ $l.\X$-${\rm Ggldim}(R)\leq n$.

$(2)$ $\X$-${\rm Gpd}(M)\leq n$ for any finitely generated $R$-module $M$.

$(3)$ $\X$-${\rm Gpd}(R/I)\leq n$ for every left ideal $I$ of $R$.

Consequently,
 \begin{align*}l.\X{\mbox -{\rm Ggldim}} (R)&={\rm sup}\{\mbox{$\X$-}{\rm Gpd}(M)\mid \mbox {$M$ is a finitely generated $R$}\mbox {-module}\}\\
&={\rm sup}\{\mbox{$\X$-}{\rm Gpd}(R/I)\mid \mbox {$I$ is a left ideal of $R$}\}.\\
\end{align*}
\end{thm}
\noindent{\bf Proof.} We need only to prove that $(3)\Rightarrow (1).$

For the case $n=0$, by the assumption, Lemma \ref{generalization of Baer criterion} shows that every $X\in\X$ is injective. It follows from \cite [Proposition 2.4] {BO} that every $R$-module is $\X$-Gorenstein
projective. Hence $l.\X$-Ggldim$(R)=0$.

Now, we can suppose that $n\geq 1$. For any $R$-module $M$, there exists a short exact sequence $0\rightarrow K\rightarrow F\rightarrow M\rightarrow0$ with $F$ free $R$-module.
We need to show that $\X$-Gpd$(K)\leq n-1$, and so $\X$-Gpd$(M)\leq n$.

 Let $F=\oplus_{i\in I}Rx_i$. We well order the index set $I$. So every element of $I$ corresponds to an ordinal which is not a limit ordinal. We add those related limit ordinals in
 $I$ to obtain a new set and still denote it by $I$. For some $k\in I$, if $k$ corresponds to a limit ordinal, we define $L_k=\oplus_{i<k}Rx_i$; if $k$ corresponds to a ordinal which
 is not a limit ordinal, define $L_k=\oplus_{i\leq k}Rx_i$. Put $J_k=K\cap L_k$, then $\X$-Gpd$(J_1)\leq n-1$ and $\X$-Gpd$(J_{i+1}/J_i)\leq n-1$ for any $1\leq i<k$. Indeed, it follows from the assumption and Corollary \ref{X-Pd of short exact sequence}, that $\X$-Gpd$(J)\leq n-1$ for every left ideal $J$ of $R$. Using the fact that $J_1=K\cap Rx_1$ and
  $J_{i+1}/J_i=K\cap Rx_{i+1}$ are isomorphic to ideals of $R$, we get that $\X$-Gpd$(J_1)\leq n-1$ and $\X$-Gpd$(J_{i+1}/J_i)\leq n-1$ for any $1\leq i<k$.
  We also note that $J_k=\cup_{i\leq k}J_i$ (resp.$J_k=\cup_{i< k}J_i$ ) for $k\in I$ which corresponds to an ordinal which is not a limit ordinal (resp. $k\in I$ which corresponds to
  a limit ordinal) and $K=\cup_{i\in I}J_i$, respectively.

  We claim that $\X$-Gpd$(J_k)\leq n-1$ for every $k\in I$.

  Firstly, we want to show that, given some $k\in I$ which corresponds to a limit ordinal, if $\X$-Gpd$(J_i)\leq n-1$ for any $i<k$, then also $\X$-Gpd$(J_k)\leq n-1$.
  From $\X$-Gpd$(J_1)\leq n-1$ and Lemma \ref{another characterization of X-Gpd}, there exists a short exact sequence
  $0\rightarrow J_1\rightarrow H_1\rightarrow G_1\rightarrow 0$, where ${\rm pd}(H_1)\leq n-1$ and $G_1$ is $\X$-Gorenstein projective.
  Consider the following pushout diagram:
 $$\xymatrix{
 & 0 \ar[d]  & 0 \ar[d] \\
 0 \ar[r] & J_1 \ar[d]^{\alpha_1} \ar[r]^{f_1} & H_1 \ar[d]^{\beta_1} \ar[r]^{g_1} & G_1 \ar@{=}[d] \ar[r] & 0  \\
  0 \ar[r] & J_2 \ar[d] \ar[r]^{l_1} & Q_2 \ar[d] \ar[r]^{h_1} & G_1  \ar[r] & 0  \\
 & J_2/J_1 \ar[d] \ar@{=}[r] & J_2/J_1 \ar[d]  \\
  & 0  & 0    }$$
  where $\X$-Gpd$(J_2/J_1)\leq n-1$ and ${\rm pd}(H_1)\leq n-1$. Then, applying Corollary \ref{X-Pd of short exact sequence} to the second exact column, we get
  that $\X$-Gpd$(Q_2)\leq n-1$, and hence $\X$-Gpd$(J_2)\leq n-1$ again by Corollary \ref{X-Pd of short exact sequence}.
  From Lemma \ref{another characterization of X-Gpd}, there exists a short exact sequence $0\rightarrow Q_2\rightarrow H_2\rightarrow T_2\rightarrow0$, where ${\rm pd}(H_2)\leq n-1$
  and $T_2$ is $\X$-Gorenstein projective. Now, using this short exact sequence and the second exact row of the diagram above, we get the following pushout diagram:
  $$\xymatrix{
       &   & 0 \ar[d]  & 0 \ar[d]   \\
        0  \ar[r] &  J_2 \ar@{=}[d] \ar[r]^{l_1} & Q_2 \ar[d]^{\beta_2} \ar[r]^{h_1} & G_1 \ar[d]^{\gamma_1} \ar[r] & 0 \\
      0  \ar[r] & J_2  \ar[r]^{f_2} & H_2 \ar[d] \ar[r]^{g_2} & H \ar[d] \ar[r] & 0  \\
       &    & T_2 \ar[d] \ar@{=}[r] & T_2 \ar[d]  \\
      & & 0  & 0   }$$
      Combining the diagrams above, we get the following exact commutative diagram:
      $$\xymatrix{
  & 0 \ar[d] & 0 \ar[d]  & 0 \ar[d]  \\
  0  \ar[r] & J_1 \ar[d]^{\alpha_1} \ar[r]^{f_1} & H_1 \ar[d]^{\theta_1=\beta_2\beta_1} \ar[r]^{g_1} & G_1 \ar[d]^{\gamma_1} \ar[r] & 0 \\
  0  \ar[r] & J_2 \ar[d] \ar[r]^{f_2} & H_2 \ar[d] \ar[r]^{g_2} & G_2 \ar[d] \ar[r] & 0 \\
  0  \ar[r] & J_2/J_1 \ar[d] \ar[r] & H_2/H_1 \ar[d] \ar[r] & G_2/G_1 \ar[d] \ar[r] & 0  \\
  & 0  & 0  & 0  }$$
 where $G_2/G_1=T_2$ is $\X$-Gorenstein projective. Note that $H_1$ and $H_2$ have finite projective dimension, so has $H_2/H_1$.
 Applying Corollary \ref{X-Pd of short exact sequence} to the bottom exact sequence of the diagram above, and using the fact that $\X$-Gpd$(J_2/J_1)\leq n-1$ and $G_2/G_1$ is $\X$-Gorenstein
 projective, we get that $\X$-Gpd$(H_2/H_1)\leq n-1$. By Proposition \ref{finite dimension equality}, we get ${\rm pd}(H_2/H_1)=\X$-Gpd$(H_2/H_1)\leq n-1$.

Recursively, we proceed to obtain the similar commutative diagram
$$\xymatrix{
  0  \ar[r]& J_i \ar[d]_{\alpha_i} \ar[r]^{f_i} & H_i \ar[d]_{\theta_i} \ar[r]^{g_i} & G_i \ar[d]_{\gamma_i} \ar[r] & 0  \\
  0 \ar[r] & J_{i+1} \ar[r]^{f_{i+1}} & H_{i+1} \ar[r]^{g_{i+1}} & G_{i+1} \ar[r] & 0   }$$
  where ${\rm pd}(H_i)\leq n-1$, ${\rm pd}(H_{i+1}/H_i)\leq n-1$, and $G_i$, $G_{i+1}/G_i$ are $\X$-Gorenstein projective for all $i+1<k$.

  Thus, we get the following short exact sequence
  $$0\rightarrow J_k=\cup_{i<k}J_i\rightarrow \cup_{i<k}H_i\rightarrow \cup_{i<k}G_i\rightarrow0.\eqno(*)$$
  By Eklof's lemma, ${\rm pd}(\cup_{i<k}H_i)\leq n-1$. We will get that $\X$-Gpd$(J_k)\leq n-1$ if we prove that $\cup_{i<k}G_i$ is a $\X$-Gorenstein projective $R$-module.
For any module $X$ in $\X$-$\mathcal{GP}(R)$, we denote $\mathbf{P}(X)$ the $\X$-complete projective resolution of $X$, that is,
$$\mathbf{P}(X):=\cdots \rightarrow\mathbf{P}(X)_{1}\rightarrow\mathbf{P}(X)_0\rightarrow\mathbf{P}(X)^0\rightarrow\mathbf{P}(X)^1\rightarrow\cdots $$
is a $\X$-complete projective resolution of $X$. Since $G_i$ and $G_{i+1}/G_i$ are $\X$-Gorenstein projective for all $i+1<k$, from the Horseshoe Lemma and the
exact sequence $0\rightarrow G_i\rightarrow G_{i+1}\rightarrow G_{i+1}/G_i\rightarrow 0 $, we can construct inductively, for each $1<i<k$, an $\X$-complete projective resolution
of $G_{i+1}$ as follows
$$\mathbf{P}(G_{i+1})=\mathbf{P}(G_{i})\oplus\mathbf{P}(G_{i+1}/G_i).$$
Then, by Eklof's lemma, the exact sequence
$$\cup_{i<k}\mathbf{P}(G_i):\cdots\rightarrow \cup_{i<k}\mathbf{P}(G_i)_1\rightarrow\cup_{i<k}\mathbf{P}(G_i)_0\rightarrow\cup_{i<k}\mathbf{P}(G_i)^0\rightarrow\cup_{i<k}\mathbf{P}(G_i)^1
\rightarrow\cdots$$
is a $\X$-complete projective resolution of $\cup_{i<k}G_i$. Hence $\cup_{i<k}G_i$ is $\X$-Gorenstein projective. Therefore, $\X$-Gpd$(J_k)\leq n-1$ by the exact sequence $(*)$.

Secondly, we prove that $\X$-Gpd$(J_k)\leq n-1$ for every $k\in I$. Indeed, if not, then the set of elements $i\in I$ such that $\X$-Gpd$(J_i)> n-1$ admits a minimal
element $k$. From the proof of the first step, $k$ can not correspond to a limit ordinal. Hence, $k-1$ exists in $I$. By the minimality of $k$, $\X$-Gpd$(J_{k-1})\leq n-1$.
 Applying Corollary \ref{X-Pd of short exact sequence}, from the fact that $\X$-Gpd$(J_k/J_{k-1})\leq n-1$, we get that $\X$-Gpd$(J_{k})\leq n-1$, which is a contradiction.
 The claim is proved.

 Finally, noting that $K=\cup_{i\in I}J_i=\cup_{i<\omega}J_i$ for some limit ordinal $\omega$, and repeating the first step to $K$, we get also that $\X$-Gpd$(K)\leq n-1$.
 We complete the proof.\hfill$\square$\\

Let $\X$ be the class of Gorenstein projective modules, then $\X$-$\mathcal{GP}(R)=\mathcal{P}(R)$ by Remark \ref{remark of some kinds of X-Gp modules}. Applying the theorem
above, we get the well-known Auslander's Theorem on global dimension.
\begin{cor}$($\cite [Theorem 1]{A55}$)$ Let $R$ be a ring. Then we have
 \begin{align*}l.{\mbox {\rm gldim}} (R)&={\rm sup}\{\mbox{\rm pd}(M)\mid \mbox {$M$ is a finitely generated $R$}\mbox {-module}\}\\
&={\rm sup}\{{\rm pd}(R/I)\mid \mbox {$I$ is a left ideal of $R$}\}.
\end{align*}
\end{cor}

Let $\X$ be the class of projective $R$-modules, then the class of $\X$-Gorenstein projective modules coincides with the class of Gorenstein projective modules. Hence we have
the following corollary.
\begin{cor}$($\cite [Theorem 1.1]{BHW}$)$ Let $R$ be a ring. Then we have
 \begin{align*}l.{\mbox {\rm Ggldim}} (R)&={\rm sup}\{\mbox{\rm Gpd}(M)\mid  \mbox {$M$ is a finitely generated $R$}\mbox {-module}\}\\
&={\rm sup}\{{\rm Gpd}(R/I)\mid \mbox {$I$ is a left ideal of $R$}\}.
\end{align*}
\end{cor}

Let $\X$ be the class of flat $R$-modules, then the class of $\X$-Gorenstein projective modules coincides with the class of strongly Gorenstein flat modules.
For an $R$-module $M$, the strongly Gorenstein flat dimension of $M$, denoted by ${\rm SGfd}(M)$, is defined as
${\rm SGfd}(M)={\rm Inf}\{n\mid \exists$ strongly Gorenstein flat resolution $0\rightarrow G_n\rightarrow \cdots G_1\rightarrow G_0\rightarrow M\rightarrow 0  \}$. We set
 ${\rm SGfd}(M)=\infty$ if there is no such a strongly Gorenstein flat resolution.

 Define $l.{\rm SGfgldim}(R)={\rm sup}\{$${\rm SGfd}(M)\mid M$ is any left $R$-module $\}$. We call it left global strongly Gorenstein flat dimension of $R$.

By Theorem \ref{X-Ggldim of R}, we have the following corollary.
\begin{cor} Let $R$ be a ring. Then we have
\begin{align*}l.{\mbox {\rm SGfgldim}} (R)&={\rm sup}\{\mbox{\rm SGfd}(M)\mid \mbox {$M$ is a finitely generated $R$}\mbox {-module}\}\\
&={\rm sup}\{{\rm SGfd}(R/I)\mid \mbox {$I$ is a left ideal of $R$}\}.
\end{align*}

\end{cor}

\noindent{\bf Acknowledgements}

The research was completed during the authors' visit at University of Washington. The authors would like to express there sincere thanks to Professor James Zhang for his hospitality.

\vspace{0.5cm}



\end{document}